\documentclass[smallextended,referee,envcountsect]{svjour3} 
\smartqed 
\usepackage{graphicx}
\usepackage{pstricks-add}
\usepackage{amsfonts,amssymb,amsmath}
\usepackage[caption=false]{subfig}
\newcommand{\rn}{\mathbb{R}^{n}}

\newcommand{\R}{\mathbb{R}}

\begin{document}

\title{Reduction of Exhausters by Set Order Relations and Cones}

%\subtitle{Reduction by Order}

\author{Mustafa Soyertem   \and  \.{I}lknur Atasever G\"{u}ven\c{c} \and Didem Tozkan }

\institute{Mustafa Soyertem Corresponding author  \at
             U\c{s}ak University \\
              U\c{s}ak, Turkey\\
              mustafa.soyertem@usak.edu.tr  
           \and
           \.{I}lknur Atasever G\"{u}ven\c{c},  \at
              Eski\c sehir Technical University \\
              Eski\c sehir, Turkey\\
              iatasever@eskisehir.edu.tr
           \and
           Didem Tozkan \at
              Eski\c sehir Technical University \\
              Eski\c sehir, Turkey\\
           dtokaslan@eskisehir.edu.tr
}

\date{Received: date / Accepted: date}
%The correct dates will be entered by the editor.

\maketitle

\begin{abstract}
The notions of upper and lower exhausters are effective tools for the study of non
smooth functions. There are many studies presenting optimality conditions for
unconstrained and constrained cases. One can observe that optimality conditions in
terms of both proper and adjoint exhausters are related to all elements of
the exhausters.

Moreover, in the constrained case the conditions that must be provided for a
particular cone determined by constraint set and the point (to be checked whether it
is optimal) are rather challenging to check. Thus it is advantageous to reduce the
number of sets in the exhauster for constrained case.

In this work, we first consider constrained optimization problems and deal with the
problem of reducing generalized exhausters of the directional derivative of the
objective function. We present some results to reduce generalized lower (upper)
exhausters by using set order relations $\preceq^{m_1}$ and $\preceq^{m_2}$, respectively. Furthermore, we
show that a generalized exhauster $E$ can be reduced to the set of minimal elements
of $E$ with respect to $\preceq^{m_1}$ or $\preceq^{m_2}$. 
 Then considering unconstrained optimization problems, lower  exhausters are reduced by using cones.
\end{abstract}

%All acknowledgements should be placed in the back of the paper after Conclusions..

\section{Introduction}
\indent\indent Exhausters defined by Demyanov \cite{d2} are geometrical and practical tools to give optimality conditions in Optimization Theory.  Exhausters and generalized exhausters are families of compact convex sets providing  representations of positively homogeneous functions and they are effectively used in nonsmooth optimization  \cite{ab2,d1,dr3,dr4}. Considering the upper and lower exhausters of the directional derivative  (in the sense of Dini, Hadamard etc.) of the objective function for a minimum or maximum problem, many optimality conditions are given \cite{d2,ab2,d1,dr2,weakexh,reduction}. Demyanov and Roshchina \cite{quasi,dr1}  gave  optimality conditions for constrained optimization problems via generalized exhausters.

%By these notions  optimality conditions are geometrically described as the directional derivatives (such as Dini, Hadamard, Clarke and Michel-Penot derivatives) which are positively homogeneous functions of directions \cite{d1}. By means of exhausters Demyanov and Roshchina \cite{dr3,dr4} presented optimality conditions for unconstrained optimization problems  and Demyanov \cite{quasi, dr1}  gave  optimality conditions for constrained optimization problems via generalized exhausters.  Furthermore, some optimality conditions in terms of lower and upper exhausters were introduced in \cite{ab2,dr1,dr2,weakexh}.

An upper/lower exhauster of a positively homogeneous function may contain
infinitely many sets. Hence researchers  have been interested in reducing the number of sets or reducing the size of the sets.  Roshchina defined minimality of exhausters both by inclusion (reducing number) and by shape (reducing the size) and gave some techniques to reduce exhausters \cite{ros2,ros3}. By using shadowing sets Grzybowski et al. \cite{gfinite} gave a criterion for minimality of finite upper exhausters. Also,  they showed that a minimal exhauster does not have to be unique.

In addition, K\"{u}\c{c}\"{u}k et al. \cite{weakexh} defined weak exhausters which are a special class of exhausters defined by weak subdifferential \cite{weaksg}. Also, they gave optimality conditions via weak exhausters \cite{weakexh} and some conditions to reduce them \cite{weakexh,reduction}. Recently, Abbasov \cite{ab1} propose new conditions for the verification of minimality of exhausters for maxmin or minmax representation of a positively homogeneous function, and present some reduction techniques.

Removing redundant sets is an effective way for reduction of an exhauster. For this purpose we use the concept of set order relations on family of sets.  Kuroiwa et al. introduced first six set order relations on family of sets  as generalizations of vector order relations \cite{Kuroiwa2}. These order relations gave a new  point of view to the concept of the solution of set-valued optimization problems. This point of view is known as ``set approach''. Furthermore, Jahn and Ha \cite{Jahn1} and Karaman et al. \cite{Karaman2} defined different set order relations. In this article, we consider $\preceq^{m_1}$ and $\preceq^{m_2}$ set orders which are defined by Karaman et al. \cite{Karaman2}. These order relations are defined by means of Minkowski difference and a cone, and they are partial order relations on the family of nonempty, bounded subsets of a normed space.

In this work, we first consider constrained optimization problems and present some reduction techniques by means of set order relations $\preceq^{m_1}$ and $\preceq^{m_2}$ for generalized exhausters of some directional derivatives of the objective function. The optimality conditions are rather challenging to check in constrained case. Thus we aim to make calculations easier by reducing the number of the sets of the generalized exhauster corresponding to a constrained optimization problem. We use set orders relation $\preceq^{m_1}$ for generalized lower exhausters and $\preceq^{m_2}$ for generalized upper exhausters, respectively. We show that the family of minimal (maximal) of a generalized lower (upper) exhauster with respect to $\preceq^{m_1}$($\preceq^{m_2}$) is still a generalized lower (upper) exhauster. Finding the family of minimal (maximal) sets with respect to $\preceq^{m_1}$($\preceq^{m_2}$)  for a given generalized lower (upper) exhauster generates a set optimization problem. So, these results enables us to consider the reduction of exhausters as a set optimization problem. We obtain these results in terms of the cone (namely $T$) and its negative dual (namely $K$) corresponding to the directional derivative which can be chosen in the sense of Dini, Hadamard, Clarke or Michel-Penot for the solution of the problem. For example, Contingent cones corresponds to Hadamard upper and lower directional derivative (see Lemma (\ref{lemma1})). Hence these results can be adapted with respect to the choice of the directional derivative and corresponding cones. Also these reduction techniques are generalized for unconstrained case. Results given via set orders in the manuscript give a new point of view of reducing exhausters.

The paper is organized as follows: In Section 2, we give some basic definitions about exhausters and set orders.  In Section 3, we present the main results on reducing generalized exhausters via set orders for constrained case. In the last section, we reduce lower exhausters for unconstrained case.

\section{Preliminaries}
\indent\indent  In this section we recall basic definitions and results used in this article.

Let $Y$ be a normed space, $D\subset Y$ be a set. $D$ is called a cone if $\lambda y\in D$ for all $\lambda\geq 0 $ and $y\in D$. The set $N(D)=\{x^*\in Y^*\ |\ x^*(y)\leq 0 \text{ for all }y\in D\}$ is called the negative dual cone of $D$ where $Y^*$ is the  dual of $Y$. Let $S\subset Y$  and $\bar x\in cl(S)$. The set 
\vspace*{-.2cm}
$$T(S,\bar x)=\{y\in Y\ |\ \exists d_n\to 0^+ \text{ and } y_n\to y\text{ such that } \bar x+d_ny_n\in S\text{ for all }n\in \mathbb N\}$$
\vspace*{-.2cm} is called the Contingent cone of $S$ at $\bar x$.

 We denote algebraic sum of the sets $A$ and $B$ by $A+B:=\{a+b \ | \ a\in A \text{ and }b\in B\}$ and Minkowski (Pontryagin) difference of $A$ and $B$ by $\displaystyle A\dot{-}B:=\{x\in X \ | \ x+B\subset A\}$.

 Throughout this paper $B(x^*,c)$ denotes the closed Euclidean ball with center $x^*$ and radius $c$, $cl(A)$ and $ co (A)$  denote the closure and convex hull of the set $A\subset \R^n$, respectively.

Now, we recall the  set order relations $\preceq^{m_1}$ and $\preceq^{m_2}$ \cite{Karaman2}. Let $D\subset Y$  be a convex, closed, pointed cone, containing $0_Y$ and have nonempty interior. We assume that $Y$ is ordered by the cone $D$.

The set order relations $\preceq^{m_1}_D$ and  $\preceq^{m_2}_D$ are given in the following definition.

\begin{definition} \cite{Karaman2}
Let $A,B$ be nonempty subsets of $Y$.
\begin{itemize}
  \item[(i)] The order relation $\preceq^{m_1}_D$ is defined by $$A\preceq^{m_1}_DB :\Longleftrightarrow (B\dot{-}A)\cap D\neq\emptyset.$$
  \item[(ii)] The order relation $\preceq^{m_2}_D$ is defined by $$A\preceq^{m_2}_DB :\Longleftrightarrow (A\dot{-}B)\cap (-D)\neq\emptyset.$$
\end{itemize}
\end{definition}

Note that $\preceq^{m_1}_D$ and  $\preceq^{m_2}_D$ are  partial order relations on  the family of nonempty and bounded subsets of $Y$ \cite{Karaman2}.

 $m_1$-minimal and $m_2$-minimal sets of a family with respect to $\preceq^{m_1}_D$ and $\preceq^{m_2}_D$ are defined as follows.

\begin{definition}\cite{Karaman2}
\label{def1}
Let $\mathcal{S}$ be a family of nonempty and bounded subsets of $Y$ and $A\in\mathcal{S}$. Then,
  \begin{itemize}
    \item[(i)] $A$ is $m_1$-minimal set of $\mathcal{S}$ if there isn't any $B\in\mathcal{S}$ such that $B\preceq^{m_1}_D A$  and $A\neq B$,
    \item[(ii)] $A$ is $m_2$-minimal set of $\mathcal{S}$ if there isn't any $B\in\mathcal{S}$ such that $B\preceq^{m_2}_D A$  and $A\neq B$,
  \end{itemize}
  The set of $m_1$-minimal and $m_2$-minimal sets of a family $\mathcal{S}$ is denoted by $m_1$-$\min \mathcal{S}$ and  $m_2$-$\min \mathcal{S}$, respectively.
\end{definition}

Now, we recall the notions of upper, lower, generalized upper and lower exhausters of a positively homogeneous function \cite{quasi}.

 Let $K\subseteq \mathbb R^n $ be a cone. A function $h:K\to\R$ is called a positively homogeneous function (p.h.) if
$h(\lambda g)=\lambda h(g)$ for all $g\in K$ and all $\lambda\geq0$.

%If $h$ is p.h. and $h(g_1+g_2)\leq h(g_1)+h(g_2)$  for all $g_1,g_2\in\rn$ i.e. $h$ is sub-additive  then $h$ is called a sublinear function. Similarly, $h$ is called a superlinear function if it is \textup{p.h.} and $h(g_1+g_2)\geq %h(g_1)+h(g_2)$  for all $g_1,g_2\in\rn$ i.e. $h$ is super-additive.   %Let $A\subseteq\rn$ be a nonempty set. The function $p_A:\rn\to\R$ defined as $$p_A(x):=\underset{a\in A}{\sup}{\langle a,x\rangle}$$ is called the \emph{support %function of the set} $A$.

Let $h:\R^n\to\R$ be a \textup{p.h. function}.
A family $E^*$ of nonempty compact convex sets in $\rn$ is called an upper exhauster of $h$ if
\begin{equation}
h(g)=\inf_{C\in E^*}\max_{v\in C}\langle v,g\rangle \text{ for all }g\in \rn.
\end{equation}
 Similarly, a family $E_*$ of nonempty compact convex sets in $\rn$ is called a lower exhauster of $h$ if
\begin{equation}
h(g)=\sup_{C\in E_*}\min_{v\in C}\langle v,g\rangle \text{ for all }g\in \rn.
\end{equation}

Let $K\subseteq \mathbb R^n $ be a cone and $h:K\to\R$ be a \textup{p.h. function}.
A family $E^*$ of nonempty closed convex sets in $\rn$ is called a generalized upper exhauster of $h$ if
\begin{equation}
h(g)=\inf_{C\in E^*}\max_{v\in C}\langle v,g\rangle \text{ for all }g\in K.
\end{equation}
 Similarly, a family $E_*$ of nonempty closed convex sets in $\rn$ is called a generalized lower exhauster of $h$ if
\begin{equation}
h(g)=\sup_{C\in E_*}\min_{v\in C}\langle v,g\rangle \text{ for all }g\in K.
\end{equation}

In order to check optimality conditions in terms of exhausters researchers prefer to deal with rather ``smaller exhausters''
. This can be obtained in two cases: We can reduce the quantity of sets or reduce the size of sets. Roshchina defined smaller exhausters by inclusion and by shape \cite{ros3}. Here we only recall definitions of smaller and minimal exhausters by inclusion needed in this work.

 \begin{definition}\cite{ros3}
    Let $h$ be a p.h. function, $E_1$ and $E_2$ be lower (upper)  exhausters of $h$. If $E_1\subset E_2$, then $E_1$ is called smaller by inclusion than $E_2$.
    \end{definition}
    \begin{definition}\cite{ros3}
    An upper (lower) exhauster $E$ of the p.h. function $h$ is called minimal by inclusion, if there is no other lower (upper) exhauster $\widetilde{E}$ which is smaller by inclusion than $E$.
    \end{definition}

For constrained optimization problems, it is enough to consider the directions in a special cone corresponding to the directional derivative to check necessary and sufficient optimality conditions. In the following definition and lemma, we recall optimality conditions given in terms of Hadamard upper, lower directional derivatives and the Contingent cone. For further optimality conditions with other directional derivatives and cones one can see \cite{quasi,bazaraa}.

\begin{definition}
  Let $f:X\to \mathbb R ,\ X\subseteq \mathbb R^n$ be an open set. Take $x\in X$ and $g\in \mathbb R^n$.
  The quantity \[f_H^\uparrow (x,g)=\limsup_{(\alpha,g') \to (0^+,g) }\dfrac{1}{\alpha}[f(x+\alpha g')-f(x)] \] is called the Hadamard upper derivative of $f$ at $x$ in the direction $g$.

  The quantity \[f_H^\downarrow (x,g)=\liminf_{(\alpha,g') \to (0^+,g) }\dfrac{1}{\alpha}[f(x+\alpha g')-f(x)] \] is called the Hadamard lower derivative of $f$ at $x$ in the direction $g$.
 \end{definition}

Let $X\subseteq \mathbb R^n$ be an open set, $f:X\to \mathbb R$ be a function and $S\subset X$ be a subset. Consider the constrained  problem $$(P)\left\{\begin{array}{l}
  \min (\max) f(x)\\
  s.t.\ x\in S.
\end{array}\right.$$
\begin{lemma}\label{lemma1}\cite{quasi}
  Let $f$ be locally Lipschitz around a point $x^*\in S$, the cone $T (S,x^*)$ be a first order uniform approximation of $S$ near the point $x^*$. If $x^*\in S $ is a local or global minimizer (maximizer) of the problem $(P)$ then
  \[f_H^\downarrow (x^*,g)\geq 0\ \  (f_H^\uparrow (x^*,g)\leq 0)\ \ \ \forall g\in T(S,x^*).\]

  If    \[f_H^\downarrow (x^*,g)>0\ \  (f_H^\uparrow (x^*,g)<0)\ \ \ \forall g\in T(S,x^*),\ g\neq 0_n\]
then $x^*$ is a strict local minimizer (maximizer) of $(P)$.

\end{lemma}
%The pair $E(h)=[E^*(h),E_*(h)]$ is called a \emph{biexhauster of the function $h$}.

Let $A\subset\rn$, $K(A)$ is the positive dual cone of $A$:
\[K(A)=\{w\in\rn:\langle w,v \rangle\geq 0, \forall v\in A\}.\]
For a cone $\Gamma$ with apex $0_n$ let 
\[\Gamma=\bigcup\{A:A\in \mathcal{A}\}\]
where $\mathcal{A}$ is a family of cones with apex $0_n$.

\begin{lemma}\label{lemma2}\cite{quasi,dr1}
Let $h:\rn\to\mathbb{R}$ be a positive homogeneous function and assume that there exists a generalized lower exhauster $E_*$ of $h$. Then the following statements are equivalent:
\begin{itemize}
\item[(i)] $h(g)\leq 0, \quad \forall g\in \Gamma$
\item[(ii)] $-(C)\cap K(A)=\emptyset,\quad \forall C\in E_*\text{ and }\forall A\in\mathcal A$
\item[(iii)] $0_n\in C+K(A),\quad \forall C\in E_*\text{ and }\forall A\in\mathcal A$
\item[(iv)] $0_n\in L_*(h,\Gamma)=\bigcap\{C+K(A):C\in E_*, A\in\mathcal A\}.$
\end{itemize}

\end{lemma}

One can see that for a problem 

\begin{equation*}
(P)\left\{
\begin{array}{rl}
\max &f(x)\\
s.t. & x\in S
\end{array}
\right.
\end{equation*}

the condition $f^\uparrow_H(x^*,g)\leq 0,\quad \forall g\in T(S,x^*)$ becomes (i) of Lemma \ref{lemma2} where $\Gamma=T(S,x^*)$ and it will be easier to check (ii), (iii) and (iv) when we reduce the generalized lower exhauster of $f^\uparrow_H(x^*,\cdot)$. 

\section{Reducing Exhausters for Constrained Case}\label{sec_red}

\indent \indent Optimality conditions in terms of directional derivatives given for constrained case are examined for the directions in a special type of cones instead of $\rn$, as mentioned above. Thus, for this case generalized upper and lower exhausters are employed as in Lemma \ref{lemma2}.

In this section, we present some reduction techniques for generalized exhausters by using $\preceq^{m_1}$ and $\preceq ^{m_2}$ order relations. In constrained problems, the cones containing the necessary directions change according to the choice of directional derivative. For example, Contingent cones are used to express optimality conditions corresponding to Hadamard upper and lower directional derivatives as in Lemma \ref{lemma1}. Thus, results of this section can be adapted to the generalized exhausters of the directional  derivatives of all senses and corresponding cones.

Throughout this paper,  $f'(\cdot,\cdot)$ denotes a directional derivative, $T$ is the corresponding cone of this directional derivative and $K$ is the negative dual cone of $T$. For instance, if one use Hadamard upper directional derivative to give optimality conditions the corresponding cone $T$  is the Contingent cone. Different directional derivatives and corresponding cones also can be used.

If any two sets of a lower exhauster are comparable with respect to $\preceq_K^{m_1}$ then the bigger one can be omitted. The following theorem shows this fact.

\begin{theorem}\label{the2}
Let $S\subset \mathbb R^n$, $f:S\to \mathbb R$, $h(g)=f'(\bar{x},g)$ for all $g\in T$, and $E_*$ be a generalized lower exhauster of $h$. If $C_1,C_2\in E_*$ and $C_1\preceq_K^{m_1} C_2$, then $\tilde{E_*}:=E_*\setminus\{C_2\}$  is a generalized lower exhauster of $h$.
\end{theorem}

Now, we give an illustrative example for Theorem \ref{the2}.
\begin{example}
\label{ex}
Let $S=\{(x_1,x_2)\in\rn:x_1\geq0, |x_2|\leq x_1\}$, $f:S\to \mathbb R$ be the function with  Hadamard upper directional derivative \linebreak $h(g_1,g_2)=|g_1|-\sqrt{g_1^2+g_2^2}$ for all $(g_1,g_2)\in\ T(S,(0,0))$. One can easily see that   $T(S,(0,0))=S$, $K=N(T(S,(0,0)))=\{(x_1,x_2)\in\rn: x_1\leq 0, |x_2|\leq -x_1 \}$. Also $E_*=\{A:=B((0,0),1),B:=B((1,0),1)\}$ is a generalized lower exhauster of $h$ (see Figure \ref{fig1} (a) and (b)). It is clear that 
\begin{eqnarray*}
\{(-1,0)\}\in (A\dot{-}B)\cap K\Rightarrow B\preceq_K^{m_1}A. \quad\text{(see Figure \ref{fig1} (c))}
\end{eqnarray*}  By Theorem \ref{the2}, we can conclude that $E_*\setminus \{A\}=\{B\}$ is also a generalized lower exhauster of $h$.
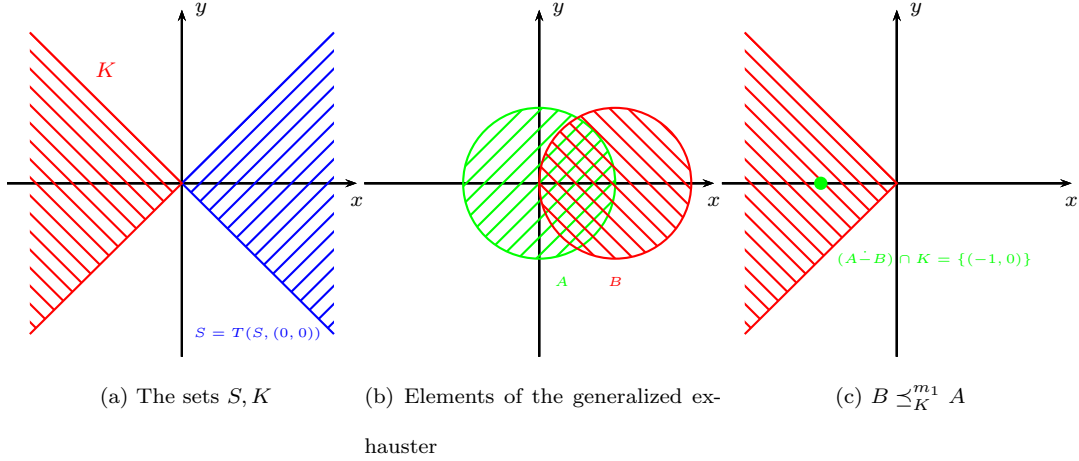
\begin{figure}[h]
     \subfloat[The sets $S,K$]{%
       \begin{pspicture}(-2.3,-2.3)(2.3,2.3)
       \psset{labels=none,ticks=none}
       \psaxes{->}(0,0)(-2.3,-2.3)(2.3,2.3)[$x$,-90][$y$,0]
       \pscustom[linestyle=none,fillstyle=hlines,hatchcolor=blue]{
		\psline(0,0)(2,-2)(2,2)(0,0)
       }
       \psline[linecolor=blue](2,-2)(0,0)(2,2)
       \pscustom[linestyle=none,fillstyle=vlines,hatchcolor=red]{
		\psline(0,0)(-2,-2)(-2,2)(0,0)
       }
       \psline[linecolor=red](-2,-2)(0,0)(-2,2)
       \rput(-1,1.5){\red $K$}
       \rput(1,-2){\tiny\blue $S=T(S,(0,0))$}
       \end{pspicture}
     }
     \subfloat[Elements of the generalized exhauster]{%
       \begin{pspicture}(-2.3,-2.3)(2.3,2.3)
       \psset{labels=none,ticks=none}
       \psaxes{->}(0,0)(-2.3,-2.3)(2.3,2.3)[$x$,-90][$y$,0]
       \pscircle[fillstyle=hlines,hatchcolor=green,linecolor=green](0,0){1}
       \pscircle[fillstyle=vlines,hatchcolor=red,linecolor=red](1,0){1}
       \rput(.3,-1.3){\tiny\green $A$}
       \rput(1,-1.3){\tiny\red $B$}
       \end{pspicture}     }
     \subfloat[$B\preceq_K^{m_1}A$]{%
       \begin{pspicture}(-2.3,-2.3)(2.3,2.3)
       \psset{labels=none,ticks=none}
       \psaxes{->}(0,0)(-2.3,-2.3)(2.3,2.3)[$x$,-90][$y$,0]
        \pscustom[linestyle=none,fillstyle=vlines,hatchcolor=red]{
		\psline(0,0)(-2,-2)(-2,2)(0,0)
       }
       \psline[linecolor=red](-2,-2)(0,0)(-2,2)

       \psdot[dotsize=5pt,linecolor=green](-1,0)
       \rput(.5,-1){\tiny \green $(A\dot{-}B)\cap K=\{(-1,0)\}$}
       \end{pspicture}     }
     \caption{Geometry of Example \ref{ex}}
     \label{fig1}
\end{figure}
\end{example}

\begin{corollary}\label{cor1} Let $S\subset \mathbb R^n$, $f:S\to \mathbb R$, $h(g)=f'(\bar{x},g)$ for all $g\in T$, and $E_*$ be a generalized lower exhauster of $h$. 
If $m_1$-$\min E_*\neq \emptyset$ then the family $m_1$-$\min E_*$ is still a generalized lower exhauster of $h$.
\end{corollary}

As a result of Theorem \ref{the2}, if a generalized lower exhauster has a strongly minimal element with respect to $\preceq^{m_1}$, then a reduced generalized exhauster consisting just this minimal element is minimal by inclusion. The following corollary states this property.

\begin{corollary}\label{the1}
Let $S\subset \mathbb R^n$, $f:S\to \mathbb R$, $h(g)=f'(\bar{x},g)$ for all $g\in T$, and $E_*$ be a generalized lower exhauster of $h$. If there exists a $C_0\in E_*$ such that \vspace*{-.3cm}\[C_0\preceq_K^{m_1} C\] \vspace*{-.3cm}for all $C\in E_*$ (i.e. $C_0$ is the strongly minimal element of $E_*$) then $\tilde{E}_*=\{C_0\}$ is a minimal exhauster by inclusion.
\end{corollary}

\begin{remark}
$ m_1$-$\min E_*$ may not be a minimal exhauster by inclusion. The following example shows this fact.
\end{remark}

\begin{example}\label{ex1}
Let $f:\mathbb R^2_+\to\mathbb R$ be a constant function and $\bar{x}=(0,0)$ where $\mathbb R^2_+$ is the first orthant. Then $T(S,\bar{x})=\mathbb R^2_+$, $K=(T(\mathbb R^2_+,\bar{x}))^*=\mathbb R^2_-$, and Hadamard upper directional derivative of $f$ at $\bar x$ is $h(g)=f^\uparrow_H(\bar x,g)=0$ for all $g\in T(\mathbb R^2_+,\bar x)$. Hence  $$E_*:=\left\{B_\alpha:=B((\cos \alpha,\sin \alpha),1):\alpha\in[0,\pi/2]\right\}$$  is a generalized lower exhauster of $h$. It is obvious that $E_*=m_1-\min E_*$. Indeed, for any $\alpha_1,\alpha_2\in[0,\pi/2]$ such that  $\alpha_1\neq\alpha_2$ 
\[B_{\alpha_1}\dot{-}B_{\alpha_2}=(\cos\alpha_1,\sin \alpha_1)-(\cos \alpha_2,\sin \alpha _2)\notin\mathbb R^2_+\] which means $B_{\alpha_1}\not\preceq_{\mathbb R^2_+}^{m_1}B_{\alpha_2}$. Since $B_{\alpha_1}$ and $B_{\alpha_2}$ are arbitrary, then we see that all elements of $E_*$ are minimal. \\
On the other hand, $\bar E_*:=\left\{B_\alpha:\alpha\in(0,\pi/2]\right\}\subsetneq E_*$ is also a generalized lower exhauster. Indeed, there exists a sequence $\{\alpha_n\}\subset(0,\pi/2]$ converges to $0$, since $0\in cl(0,\pi/2]$. Then $\displaystyle \lim_{n\to\infty}\underset{v\in B_{\alpha_n}}{\min}\langle v,g\rangle=\underset{v\in B_{0}}{\min}\langle v,g\rangle$ for all $g\in\mathbb R^2_+$. Hence,
\[\underset{\alpha\in (0,\pi/2]}{\sup}\underset{v\in B_\alpha}{\min}\langle v,g\rangle=\underset{\alpha\in [0,\pi/2]}{\sup}\underset{v\in B_\alpha}{\min}\langle v,g\rangle.\] Hence $B_0=B((1,0),1)$ can be omitted from $E_*=m_1$-$\min E_*$. Thus, $m_1$-$\min E_*$ is not minimal by inclusion.
\end{example}
As $$A\preceq_K^{m_1}B \iff -B\preceq_K^{m_2}-A  $$ similar reducing results can be obtained for generalized  upper exhauster by using $\preceq_K^{m_2}$ as follows. 

\begin{theorem}\label{the4}
Let the assumptions of Theorem \ref{the2} are valid and $E^*$ be a generalized  upper exhauster of $h$.
 If $C_1,C_2\in E^*$ satisfies $C_1\preceq_K^{m_2} C_2$ then $\bar{E}^*=E^*\setminus\{C_1\}$ is also a generalized upper exhauster of $h$.
\end{theorem}
\begin{corollary}\label{cor2}
Let $S\subset \mathbb R^n$, $f:S\to \mathbb R$, $h(g)=f'(\bar{x},g)$ for all $g\in T$, and $E^*$ be a generalized upper exhauster of $h$. 
If $m_2$-$\max E^*\neq \emptyset$ then the family $m_2$-$\max E^*$ is still a generalized upper exhauster of $h$.
\end{corollary}
\begin{corollary}\label{the3}
 Let the assumptions of Theorem \ref{the4} are valid. If there exists a set $C_0\in E^*$ such that $C\preceq_K^{m_2} C_0$ for all $C\in E^*$, then $\bar{E}^*=\{C_0\}$ is a generalized upper exhauster of $h$.
\end{corollary}

\section{Reducing Exhausters for Unconstrained Case}
All the results given in previous section are presented for generalized exhausters for constrained case. Here  we consider unconstrained optimization problems where $\rn$ becomes corresponding cone $T$. The results given in Section 3 can be generalized for lower exhausters as follows.

\begin{theorem}\label{the5}
Let $f:\mathbb R^n\to \mathbb R$, $\bar{x}\in\mathbb R^n$, $h(g):=f'(\bar{x},g)$ for all $g\in\mathbb R^n$, $E_*$ be a lower exhauster of $h$ and  $A\in E_*$. If there exist $B_1,B_2,\cdots,B_m\in E_*$ satisfying \[\bigcup_{i=1}^m(A\dot{-}B_i)^\#=\mathbb R^n\] where $\displaystyle (A\dot{-}B_i)^\#:=\!\!\!\bigcup_{d\in A\dot{-}B_i}\{x\in\mathbb R^n:\langle x,d\rangle\leq 0\}$, then $E_*\setminus\{A\}$ is a lower exhauster of $h$.
\end{theorem}

Now we give an illustrative example for this theorem.

\begin{example}
\label{ex2}
Let $h:\mathbb R^2\to \mathbb R$ be defined as $h(g_1,g_2)=|g_1|-\sqrt{g_1^2+g_2^2}$ for all $(g_1,g_2)\in\mathbb R^2$. Then $E_*=\{A:=B((0,0),1),B_1:=B((-1,0),1),B_2:=B((1,0),1)\}$ is a lower exhauster of $h$ (see Figure \ref{fig2}(a)). We get
\begin{eqnarray*}
(A\dot{-}B_1)^\#&=&\{(1,0)\}^\#=\{(x,y)\in\mathbb R^2:x\leq 0\}\\
(A\dot{-}B_2)^\#&=&\{(-1,0)\}^\#=\{(x,y)\in\mathbb R^2:x\geq 0\}.
\end{eqnarray*}

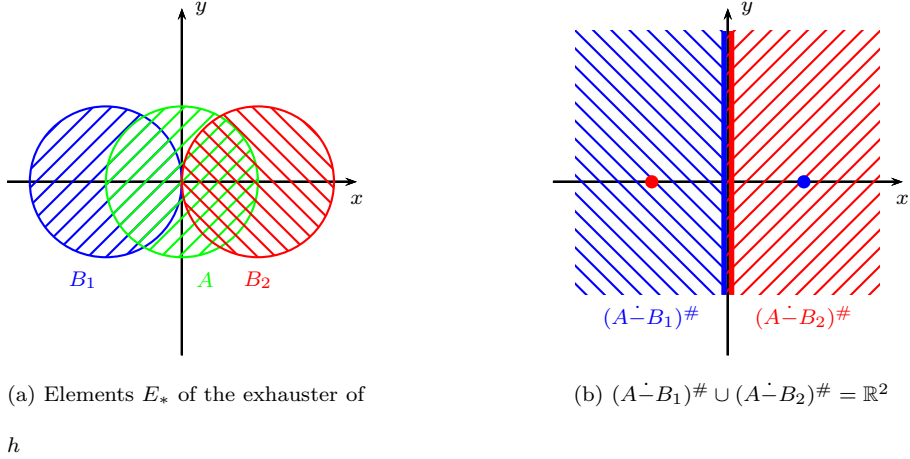
\begin{figure}[h]
     \subfloat[Elements $E_*$ of the exhauster of $h$]{%
       \begin{pspicture}(-2.3,-2.3)(2.3,2.3)
       \psset{labels=none,ticks=none}
       \psaxes{->}(0,0)(-2.3,-2.3)(2.3,2.3)[$x$,-90][$y$,0]
       \pscircle[fillstyle=hlines,hatchcolor=blue,linecolor=blue](-1,0){1}
       \pscircle[fillstyle=hlines,hatchcolor=green,linecolor=green](0,0){1}
       \pscircle[fillstyle=vlines,hatchcolor=red,linecolor=red](1,0){1}
       \rput(-1.3,-1.3){\blue $B_1$}
       \rput(.3,-1.3){\green $A$}
       \rput(1,-1.3){\red $B_2$}
       \end{pspicture}     }
       \hfill
     \subfloat[$(A\dot{-}B_1)^\#\cup(A\dot{-}B_2)^\#=\mathbb R^2$]{%
       \begin{pspicture}(-2.3,-2.3)(2.3,2.3)
       \psset{labels=none,ticks=none}
       \psaxes{->}(0,0)(-2.3,-2.3)(2.3,2.3)[$x$,-90][$y$,0]
        \pscustom[linestyle=none,fillstyle=vlines,hatchcolor=blue]{
		\psline(0,-1.5)(0,2)(-2,2)(-2,-1.5)(0,-1.5)
       }
       \psline[linecolor=blue,linewidth=2pt](-.05,-1.5)(-.05,2)
        \pscustom[linestyle=none,fillstyle=hlines,hatchcolor=red]{
		\psline(0,-1.5)(0,2)(2,2)(2,-1.5)(0,-1.5)
       }
       \psline[linecolor=red,linewidth=2pt](.05,-1.5)(.05,2)
       \psdot[dotsize=5pt,linecolor=red](-1,0)
       \psdot[dotsize=5pt,linecolor=blue](1,0)
       \rput(1,-1.8){\red $(A\dot{-}B_2)^\#$}
       \rput(-1,-1.8){\blue $(A\dot{-}B_1)^\#$}
       \end{pspicture}     }
     \caption{Geometry of Example \ref{ex2}}
     \label{fig2}
\end{figure}

It is clear that $(A\dot{-}B_1)^\#\cup(A\dot{-}B_2)^\#=\mathbb R^2$ (see Figure \ref{fig2}(b)) . By Theorem \ref{the5}, we can conclude that $E_*\setminus \{A\}=\{B_1,B_2\}$ is also a lower exhauster of $h$.
\end{example}

In Example \ref{ex2}, we see that  $M_*(A)=clco\{(0,-1),(0,1)\}\neq\emptyset$ that means this lower exhauster $E_*$ cannot be reduced to $E_*\setminus \{A\}$ by Theorem 2.6 in \cite{ros3}. Altough the method in this manuscript requires the nonemptiness of Minkowski difference of some of sets in the generalized exhausters, this example shows that some generalized exhausters which can not be reduced via other methods in the literature can be reduced via this method.

\begin{lemma}\label{lem1}
Let $A,B\subset \mathbb R^n$ compact, convex sets. The following statements are equivalent:
\begin{itemize}
\item[(i)] $(A\dot{-}B)^\#=\mathbb R^n$
\item[(ii)] $0\in A\dot{-}B$
\item[(iii)] $B\subset A$
\end{itemize}
\end{lemma}

\begin{remark}
   If $(A\dot{-}B)^\#=\mathbb R^n$ for $A,B\in E_*$, then from  Lemma \ref{lem1} we have $B\subset A$. Therefore, by Theorem \ref{the5} the same result with  Theorem 4.4 (i) in \cite{ros2} is obtained.
\end{remark}

\end{document}